\documentclass{article}
\usepackage{amsmath,amsthm,amsfonts,graphicx}
\newtheorem{theorem}{Theorem}
 
\theoremstyle{definition}
 
\newtheorem{algorithm}[theorem]{Algorithm}

\newtheorem{remark}[theorem]{Remark}

\def\smallddots{\mathinner{\raise7pt\hbox{.}\raise4pt\hbox{.}\raise1pt\hbox{.}}} 
\def\smallsdots{\mathinner{\raise1pt\hbox{.}\raise4pt\hbox{.}\raise7pt\hbox{.}}}

\begin{document}
\title{Root-finding with Implicit Deflation}
\author{Victor Y. Pan} 
\author{R\'emi Imbach$^{[1],[a]}$,
Victor Y. Pan
 $^{[2, 3],[b]}$,
Chee Yap$^{[1],[c]}$, \\
Ilias S. Kotsireas$^{[4],[d]}$ and 
Vitaly Zaderman$^{[3],[e]}$
\\
\and\\
$^{[1]}$ 
Courant Institute of Mathematical Sciences\\
 New York University,USA\\
$^{[2]}$
Department of Computer Science \\
Lehman College of the City University of New York \\
Bronx, NY 10468 USA and \\
$^{[3]}$ 
Ph.D. Programs in Mathematics  and Computer Science \\
The Graduate Center of the City University of New York \\
New York, NY 10036 USA \\
$^{[4]}$
Wilfrid Laurier University \\
Department of Physics  
and Computer Science  \\
75 University Avenue West  \\
Waterloo, Ontario N2L 3C5
CANADA  \\
$^{[a]}$
remi.imbach@nyu.edu \\
$^{[b]}$
 victor.pan@lehman.cuny.edu \\
http://comet.lehman.cuny.edu/vpan/ \\
$^{[c]}$  Email: yap@cs.nyu.edu www.cs.nyu.edu/yap/ \\
$^{[d]}$ ikotsire@wlu.ca \\
http://web.wlu.ca/science/physcomp/ikotsireas/\\
$^{[e]}$  vza52@aol.com \\} 

\date{} 
\maketitle  
 

\begin{abstract}   
Functional iterations such as Newton's are a popular tool for polynomial root-finding. We consider a realistic situation where some roots have already been approximated (we say tamed), and one would like to restrict further root-finding to the approximation of the remaining (wild) roots. 
A natural approach of applying 
explicit deflation has been much studied and recently advanced by one of the authors of this paper, but presently we consider the alternative of implicit deflation combined with  mapping of the variable and reversion of an input polynomial. The hope is that the union of the sets
of tame roots approximated in a number of such transformations can cover all roots of a polynomial.
 
 We also show another direction to substantial further progress in this long and extensively studied area. Namely we dramatically increase the local efficiency
 of root-finding  by means of the incorporation of fast algorithms for multipoint polynomial evaluation and the Fast Multipole Method.
\end{abstract} 






\section{Introduction}\label{s1}


 Univariate polynomial root-finding, that is, approximation of the roots 
$x_{1}$, $\dots$, $x_d$  of a polynomial equation
\begin{equation}\label{eqpolyp}
p(x)=0~{\rm for}~ 
p(x)=\sum^{d}_{j=0}p_jx^j=p_d\prod^d_{i=1}(x-x_i),~~~ p_d\ne 0,
\end{equation}  
has been the central problem of Mathematics for four millennia, since  the Sumerian times.
It is still involved in various areas of modern computation and is the subject of intensive research worldwide. The user's choice since 2000 has been the package MPSolve (cf. \cite{BF00}, \cite{BR14}), which implements
Ehrlich's functional iterations, but other
functional iterations such as Newton's and Weierstrass's are also highly popular.
Ehrlich's and Weierstrass's iterations 
converge simultaneously to all complex roots of a polynomial.  Newton's iterations
converge to  a single root but can be extended to approximation of all roots
or the roots in a fixed domain.

Usually root-finding iterations  approximate (we say {\em tame}) most of the roots, and then one can deflate an input polynomial and keep updating only the approximations to the remaining, {\em wild} roots.
  
Efficient methods for explicit deflation can be found in   \cite{P19} and references therein, but here we  study  alternative  techniques of implicit deflation, which enable us to exploit the sparseness of an input and to avoid numerical stability problems caused by the coefficient growth in factorization of a polynomial. 

The partition of the root set into tame and wild roots for a fixed root-finder depends on the roots disposition relatively to the initial root approximations. So the partition should change if we change
the initial approximations or if we map the variable, apply the same root-finder to the new resulting polynomial, and then recover the roots of the original polynomial by  applying the converse map.
 By using implicit deflation we avoid computing the new coefficients in these mappings
 (see Sections \ref{s3} and \ref{snnlrmp}).
We  conjecture that already the union of a small number of the resulting variations of the set of tame roots would include all roots of $p$, and we can strengthen our chances for success of this heuristic approach by  applying it  concurrently using various  functional iterations.

 In Section \ref{seff} we point out another promising direction to enhancing the power of root-finding iterations, namely by means of incorporation of superfast multipoint polynomial evaluation and the Fast Multipole Method. We demonstrate the high promise of this approach by showing that it yields a dramatic increase of local efficiency of root-finding iterations. 

Otherwise we organize our paper as follows.
In the next section we recall some popular functional iterations for  polynomial root-finding. In Section \ref{sextw}
we comment on   partitioning polynomial roots into
tame ones (already approximated) and wild ones. In Section \ref{s2} we compare explicit and implicit deflation and specify  implicit deflation for Newton's iterations.
 We combine implicit deflation  with linear maps of the variable and reversion of a polynomial in Section \ref{s3} and with squaring the variable in Section \ref{snnlrmp}, followed by our  comments on potential benefits of 
 concurrent  root-finding   in Section \ref{stmw}. 


\section{Functional Iterations for Root-Finding}\label{seffrf}


Among hundreds if not thousands known  polynomial root-finders (see up to date coverage in 
\cite{M07},  \cite{MP13}, \cite{P19}, and the bibliography therein) consider the class of functional iterations. For a fixed set of functions
$$f_1(z),\dots,f_m(x),~1\le m\le d,$$
these iterations recursively refine  current approximations $z_1^{(k)},\dots,                                                                                                                                                                                                                                                                        z_m^{(k)}$
to $m$ roots $x_1,\dots x_m$ of $p(x)$ according to the expressions 
\begin{equation}\label{eqfnit}
z_i\leftarrow f_i(z_i),~i=1,\dots,m. 
\end{equation} 
 In the case where $m=1$
  write $f(z)=f_1(z)$
 and 
 \begin{equation}\label{eqfnit1}
z\leftarrow f(z). 
\end{equation}  
These iterations
include 
 various interpolation methods, which use no derivatives of $p(x)$ and are recalled in \cite[Section 7]{MP13}, for example, Muller's method
 (see \cite[Section 7.4]{MP13});
 methods involving derivative
   such as  Newton's iterations
   \cite[Section 5]{M07}; 
and methods involving  
  higher order derivatives   
   \cite[Section 7]{MP13}.  
 We exemplify our study  
 with  Newton's iterations (where $m=1$):
 \begin{equation}\label{eqnewt}
z\leftarrow z-N_p(z),
\end{equation}
 \begin{equation}\label{eqnewt1}
~N_p(x)=p(x)/p'(x),
\end{equation}
which have efficient extensions to
the solution of polynomial
systems of equations \cite{BCSS98}
and to root-finding for various 
smooth functional equations and systems of
equations \cite{E94}; 
 Weierstrass's iterations of \cite{W03}
 (rediscovered
by Durand in \cite{D60} and  Kerner in \cite{K66}), in which case $m=d$ :
 \begin{equation}\label{eqwdk}
z_i\leftarrow z_i-W_{p,l}(z_i),~i=1,\dots,d,
\end{equation}
 \begin{equation}\label{eqwdk1}
 W_{p,l}(x)=\frac{p(x)}{p_n l'(x)}, 
 \end{equation}
  \begin{equation}\label{eql}
 l(x)=\prod_{i=1}^d(x-z_i), 
\end{equation}
and  Ehrlich's iterations of \cite{E67}
(rediscovered by Aberth in \cite{A73}),
where again $m=d$:
\begin{equation}\label{eqehrab}
z_i\leftarrow z_i-E_{p,i}(z_i),
\end{equation}
\begin{equation}\label{eqe}
E_{p,i}(x)=0~{\rm if}~p(x)=0;~
\frac{1}{E_{p,i}(x)}=\frac{1}{N_p(x)}-
\sum_{j=1,j\neq i}^d\frac{1}{x-z_j} 
~{\rm otherwise};
\end{equation}
$i=1,\dots,d$, and $N_p(x)$ is defined by (\ref{eqnewt1}). 
\begin{remark}\label{resec}
The above root-finders are readily  extended to any function $s(x)$ sharing
its  root set with the polynomial $p(x)$.
For example, deduce from  the Lagrange interpolation formula that
 $$p(x)=l(x)s(x),$$
 $$s(x)=p_n+\sum_{i=1}^d \frac{W_{p,l}(z_i)}{x-z_i}$$
  for any set of $d$ distinct
  nodes $z_1,\dots,z_d$.
  Apply selected iterations to the above {\em secular rational  function} $s(x)$ or the polynomial $l(x)s(x)$.
Bini and Robol in \cite{BR14} show  substantial benefits of that application of Ehrlich's iterations to
$l(x)s(x)$ rather than $p(x)$, both for convergence acceleration and error estimation.
\end{remark} 


\section{Tame and Wild Roots}\label{sextw}


Now suppose that we have applied a fixed functional iteration (\ref{eqfnit}) and have approximated 
$m$ roots of a polynomial $p(x)$ for $m<d$
(we call them {\em tame}); 
next we discuss efficient approximation of the  
 remaining  roots; we call them {\em wild} and call their approximation {\em taming}. 
 
For example, we face a taming problem where  functional
iterations (\ref{eqfnit1}) have approximated a single root of a polynomial $p(x)$ and we seek the other roots. 

Newton's and other iterations
(\ref{eqfnit1}), devised for approximation of a single root, can be also applied at a number of initial points in order to approximate all roots. This can succeed for most of the roots, while some  roots can escape and stay wild. In particular in  the paper  \cite{SS17}  Newton's iterations initialized at a universal
set of $O(d)$ points\footnote{This set is {\em universal} for all polynomials $p(x)$ 
that have all roots lying in the unit disc $D(0,1)=\{z:~|z|=1\}$. Given any polynomial $p(x)$ one can move all its roots into this disc by means of first readily computing a reasonably close upper bound on the absolute values of all roots and then properly shifting and scaling the variable $x$.} approximate
 $t=d-w$ roots of $p(x)$ but leave out a narrow set 
 of $w$ wild roots  where
  $w<0.001~d$ for $d<2^{17}$ and $w<0.01~d$ for $d<2^{20}$. The paper   \cite{SS17}    
 continued a long study traced back to   \cite{KS94} and \cite{HSS01}.

Likewise some roots remain wild 
while most of the roots are tamed in Weierstrass's,
 Ehrlich's, and various other iterations that recursively  update approximations of all roots.

Finally the subdivision  root-finding
iterations  of \cite{BSSY18} 
  extend the earlier study in \cite{W24}, \cite{HG69},
 \cite{H74}, \cite{R87}, and \cite{P00},
 where such iterations are called  the Quad-tree construction. This  root-finder has recently been implemented in  \cite{IPY18}. It  first approximates some sets of tame roots of $p(x)$ in certain  domains on the complex plane well-isolated from the other  roots and then approximates the remaining wild roots, in particular by combining the subdivision process with complex extension of Abbott's real QIR iterations.
   
 
\section{Taming Wild Roots by Means of  Deflation}\label{s2}


Seeking wild  roots one can deflate an input polynomial, that is, apply a selected root-finder 
to the polynomial
\begin{equation}\label{eqpolyq}
 q(x)=\sum^{w}_{i=0}q_ix^i=p_d\prod^w_{j=1}(x-x_j),~~~ p_d\ne 0.
\end{equation} 

In  {\em explicit deflation} we first compute the coefficients of $q(x)$.
If the roots of the quotient 
$q(x)$ are well isolated from the other
roots of $p(x)$, we can apply the
 efficient method of Delves and Lyness \cite{DL67}. The root-finders of \cite{S82} and \cite{K98} incorporate its advanced versions; \cite{P19} presents them in a concise form.

Bini and Fiorentino argue in \cite{BF00} that explicit deflation of a polynomial $p(x)$ does not preserve its sparseness and in some cases can be
  numerically unstable, for instance, in the case of a polynomial
  $p(x)=x^d\pm 1$ of a large degree $d$.
  These potential problems somewhat limit the value of explicit deflation  where a polynomial $q(x)$ has large degree $w$.
 Moreover we can completely avoid these problems 
by applying {\em implicit
deflation}, that is, applying functional iterations that evaluate  $q(x)$ at a point $x$  as the ratio $\frac{p(x)}{t(x)}$ for $t(x)=p_d\prod^d_{j=1+w}(x-x_j)$. 

We can readily implement this recipe in the case of functional interpolation iterations of \cite[Section 7]{MP13}. 

Let us specify implicit deflation when we apply Newton's iterations and the following well-known identity (cf. \cite{M54}),
\begin{equation}\label{eqnrt1}
 \frac{1}{N_p(x)}=\sum^n_{j=1}\frac{1}{x-x_j}.
\end{equation} 


\begin{algorithm}\label{alg1} {\em Implicit Deflation with Newton's iterations.}


\begin{description}


\item[{\sc Input:}] 
A polynomial $p(x)$ of (\ref{eqpolyp}),
 a set of sufficiently close approximations\footnote{We assume that we can very quickly refine approximations to tame roots.} to its tame roots $x_{w+1},\dots,x_d$,
an initial approximation $z$ to a wild root of $p(x)$, a Stopping Criterion (see, e.g., \cite{BF00}, \cite{BR14}), 
and a black-box program ${\rm EVAL}_p$ that evaluates the ratio $\frac{1}{N_p(z)}=\frac{p'(z)}{p(z)}$ for
 a polynomial $p(x)$ of (\ref{eqpolyp}) 
 and a complex point $z$.


\item[{\sc Output:}] 
The updated approximation  $z\leftarrow z-N_p(z)$ to a root of $p(x)$ (see (\ref{eqnewt})). 




\item[{\sc Computations}:]
Apply Newton's iteration (\ref{eqnewt}) to the polynomial $q(x)$ defined implicitly, that is, successively compute
the values:
\begin{enumerate}
\item 
$r=p'(z)/p(z)\leftarrow 1/N_p(z)$,
\item 
 $s\leftarrow \sum_{j=w+1}^d \frac{1}{z-x_j}$,
\item 
  $N_q(z)=\frac{q(z)}{q'(z)}\leftarrow \frac{1}{r-s}$.
\item 
$z\leftarrow z-N_p(z)$.
\item 
If the fixed Stopping Criterion is met, output $z$ and stop.
Otherwise  go to stage 1.
\end{enumerate}
\end{description}


\end{algorithm}


Dario A. Bini (private communication) proposed to 
improve numerical stability of this algorithm by means of scaling as follows:  

$$N_q(z)=\frac{1/r}{1-s/r}.$$

\medskip 

{\bf Complexity of a single iteration of Algorithm \ref{alg1}.}

 Stage 1 amounts to a single  invocation of the program ${\rm EVAL}_{p}$.

Stage 2 involves
$d-w$ divisions, $d-w$,
subtractions
and $d-w-1$ additions.
   
Stages 3 and 4 together involve
 $2$ subtractions  and  a single division. 
 
We can readily extend implicit deflation to various other root-finders involving Newton's ratio $N_p(x)$, in particular,
to Ehrlich's iterations (\ref{eqehrab}) because we can assume that
$z_j=x_j$ for $j>w$, and then
equation (\ref{eqnrt1})  implies that
$E_{p,j}(x)=E_{q,j}(x)$ for
$q(x)$ of (\ref{eqpolyq}), $E_{p,j}(x)$ of 
(\ref{eqe}),  and $j\le w$.


\section{Combining Newton's Iterations, Linear Maps 
of the Variable and Reversion}\label{s3}

 
The set of tame roots output 
by  fixed functional iterations  
 varies when an input polynomial $p(x)$  varies. This suggests that we can 
approximate many or all  wild roots if we
reapply the same iterations to the
polynomials
\begin{equation}\label{eqvabc}
v(z)=v_{a,b,c}(z)=(z+c)^d p\Big (a+\frac{b}{z+c}\Big )
\end{equation}
 for  
various  triples of complex scalars 
$a$, $b\neq 0$, and $c$.
We must limit the overall number of the triples  
in order to control the overall computational cost.

The following equations map the roots $x_j$ of $p(x)$
to the roots 
$z_j$ of $v(x)$ and vice versa,
\begin{equation}\label{eqrrr}
x_j=a+\frac{b}{z_j+c},~z_j=\frac{b}{x_j-a}-c.
\end{equation}
 
Let us specify this recipe
for the algorithm of \cite{SS17},  cited in Section \ref{sextw}.

\begin{algorithm}\label{algmap}
\item[{\sc Initialization}:]  
Define a polynomial $v(z)=v_{a,b,c}(z)$ by
choosing the parameters $a$, $b$, and $c$
such that all roots of the polynomial $v(z)$ lie
in the unit disc $D(0,1)=\{z:~|z|=1\}$;
 do not actually compute the coefficients of that polynomial.
\item[{\sc Computations}:]
\begin{enumerate} 
\item
Apply Newton's iteration (\ref{eqnewt}) to the  polynomial $v(z)$
by using initialization at the uni\-ver\-sal set
of \cite{SS17} and 
 by expressing
 the Newton's ratios $N_v(z)=v(z)/v'(z)$ 
(cf. (\ref{eqnewt})) via the following equations:
\begin{equation}\label{eqrmp}
\frac{1}{N_v(z)}=\frac{d}{z+c}-\frac{b}{(z+c)^2N(x)}~{\rm for}~
v(z)~{\rm of~(\ref{eqvabc})~and~}x~{\rm of~(\ref{eqrrr})}.
\end{equation}
\item
Having approximated a root $z_j$ of $v(z)$ for any $j$,
 readily recover the root $x_j$ of $p(x)$ from equation
(\ref{eqrrr}).
\end{enumerate}
\end{algorithm}

In the particular case where $a=c=0$ and $b=1$, 
the above expressions are simplified:
$z=1/x$; $v(z)$ turns into the reverse  polynomial of $p(x)$,
$$v(z)= p_{\rm rev}(z)=\sum^{d}_{i=0}p_{d-i}z^i=z^dp(1/z),$$
 $$\frac{1}{N_v(z)}=\frac{v'(z)}{v(z)}=\frac{d}{z}-\frac{1}{z^2N_p(1/z)},$$
 and $p_{\rm rev}(x)=p_0\prod^d_{j=1}(x-1/x_j)$ if $p_0\neq 0$. 


\section{Combining Newton's Iterations and Squaring of the Variable}\label{snnlrmp}


One can hope to obtain all roots of $p(x)$
by applying Newton's iterations to the polynomials $v(z)=v_{a,b,c}(z)$
 for a reasonable number of triples 
of $a$, $b$, and $c$,
but one can also extend this
approach by using
more general rational maps $y=r(x)$  (cf., e.g., \cite{MP00}).

For a simple example, 
consider the Dandelin's root-squaring
map of 1826, rediscovered by 
Lobachevsky in 1834 and then by 
Gr{\"a}ffe in 1837 (see \cite{H59}):   

\begin{equation}\label{eqdnd}
u(y)=(-1)^dp(\sqrt{y})p(-\sqrt{y})=\prod_{j=1}^d (y-x_j^2).
\end{equation} 
In this case one should  make a polynomial $p(x)$ of (\ref{eqpolyp})  monic by scaling the variable $x$ and 
then express the Newton's ratio 
$N_u(y)=u(y)/u'(y)$ as follows:

$$\frac{1}{N_u(y)}=
0.5 \Big (\frac{1}{N_p(\sqrt y)}-\frac{1}{N_p(\sqrt {-y})}\Big )~y^{-1/2}.$$
Notice that under map (\ref{eqdnd}) the roots lying in the unit disc $D(0,1)$ stay in it.

Having approximated the $n$ roots $y_1,\dots,y_n$ of the polynomial $u(y)$,
we readily recover    
 the $n$ roots $x_1,\dots,x_n$ of the polynomial $p(x)$ by selecting them
from the $2n$ values $\pm \sqrt {y_j}$, $j=1,\dots,n$. 

We can 
 apply the above maps recursively (a limited number of times,
 in order to control the overall computational cost);
 then we can recover the roots from their images in these rational maps
by extending the lifting/descending techniques of \cite{P95}, \cite{P02}.


\section{Concurrent Root-finding}\label{stmw}


\begin{remark}\label{reprl}
We recalled that Newton's iterations can compute
most of the roots of a fixed polynomial but not all of them.
Seeking the remaining, wild roots, we applied the iterations to a number of related polynomial. This recipe can be immediately extended to  application of Ehrlich's, Weierstrass's or another fixed iterative  root-finder to a variety of polynomials linked to 
an input polynomial. 
Furthermore we can 
extend this idea to
concurrent application of
a number of iterative root-finders to such a variety of  polynomials, and
 one can perform  computations on a number of processors  with minimal need for 
their communication and synchronization.
\end{remark}

\begin{remark}\label{repr2}
Weierstrass's and Ehrlich's
  functional iterations,
as well as their Gauss-Seidel's and Werner's
 accelerated  variations (cf. \cite{BR14} and \cite{W82}) 
converge very fast empirically, but formal support of this empirical observation is a well-known challenge. Can we facilitate obtaining such a support if we allow random maps of the  variable $x$, e.g., if
we apply these iterations to the polynomials $v_{a,b,c}(z)$ of (\ref{eqvabc}) for random choice of the parameters $a$, $b$, and $c$?  For example, initialization of Newton's iterations at a set of points 
$\{c+r\exp(\phi_j{\bf i})$, $j=1,\dots,s$,  
of a  circle $\{x:~|x-c|=r\}$ on the complex plane can be equivalently interpreted as  the application of these iterations at a single point $y=c$  to a set of polynomials 
$p_j(y)$ obtained from $p(x)$ via  the linear maps $y\leftarrow x-r\exp(\phi_j{\bf i})$,
$j=1,\dots,s$.  
\end{remark}


\section{Efficiency of Root-finding Iterations}\label{seff}


Since Ostrowski's paper \cite{O66}, it is customary  to measure local
efficiency of functional root-finding iterations by the quantity eff$=
q^{1/\alpha}$
or sometimes  
$\log_{10}({\rm eff})=(1/\alpha)\log_{10}q$ where $q$ denotes the convergence order (rate) and $\alpha$ is the number of function evaluations per iteration and per root. In particular $q=2$, 
$\alpha=2$,  and 
eff=$\sqrt 2\approx 1.414$ for 
Newton's and  Weierstrass's iterations while $q=3$, 
$\alpha=3$,  and eff=$3^{1/3}\approx 1.442$ for  Ehrlich's iterations where we assign the same cost to the evaluation of the functions $\sum_{j=1,j\neq i}^d\frac{1}{x-z_j}$, $p(x)$,  $p'(x)$,
 and $l'(x)$ at $x=z_i$,
 noting that $l'(z_i)=
 \prod_{j=1,j\neq i}^d(z_i-z_j)$. 

Actually the cost of function evaluation requires further elaboration.
 Exact evaluation  of  the values
$\sum_{i=1,i\neq j}^d\frac{1}{z_j^{(k)}-z_i^{(k)}}$ for $j=1,\dots,d$ 
is Trummer's  celebrated problem, whose solution, 
like exact evaluation  of a polynomial $p(x)$ of (\ref{eqpolyp}) at 
$d$ points, 
involves $O(d\log^2(d))$ arithmetic
operations \cite[Section 3.1]{P01}, \cite{GGS87}, \cite{MB72}.

Both of these superfast algorithms 
-- for  polynomial evaluation and the 
Trummer's problem --
are numerically unstable for $d>50$, but one can use numerically stable superfast alternatives  based on the Fast Multipole  Method \cite{BY13}.
Its application to Trummer's problem
 is well-known \cite{GR87}, but in the case of multipoint polynomial evaluation is more recent  and more                                                                                                                                                                                       involved    \cite{P15}  and \cite{P17}.
 
Using superfast algorithms for both problems  decreases
$\alpha$  to the  order of $O(\log^2(d)/d)$. Hence local efficiency of  Weierstrass's and  Ehrlich's iterations grows to   infinity as $d\rightarrow \infty$, and similarly for Newton's  iterations initialized  and applied  simultaneously at the 
 order of $d$ points.

The above formal analysis applies locally,
where the convergence to the roots becomes superlinear, while  the overall computational cost is usually dominant at the previous initial stage, for which only limited formal results are available
(see also  Remark \ref{repr2}). 
These limited results favor Ehrlich's
iterations, which empirically have
   milder
 sufficient conditions for superlinear
convergence  
than both Newton's and Weierstrass's iterations \cite{T98}.


 \medskip 

\noindent {\bf Acknowledgements:}
The research of R. Imbach,
V. Y. Pan, V. Zaderman and C. Yap  was supported by NSF Grant CCF-1563942. The research of R. Imbach was also supported by NSF Grants  CCF-1564132 and CCF-1708884. The research of V. Y. Pan and V. Zaderman was also supported by NSF Grant  CCF 1116736  and PSC CUNY Award 69813 00 48. The research of Ilias Kotsireas was supported by an NSERC grant.
  

\end{document}